\documentclass[11pt,a4paper]{article}

\usepackage{amsmath,amsfonts,amssymb,latexsym,graphics,epsfig,url}
\usepackage{color}
\usepackage{amsthm}
\usepackage[english]{babel}
\usepackage{graphicx}
\usepackage{soul}
\usepackage[utf8]{inputenc}

\newtheorem{theorem}{Theorem}[section]
\newtheorem{proposition}[theorem]{Proposition}
\newtheorem{lemma}[theorem]{Lemma}
\newtheorem{corollary}[theorem]{Corollary}

\newtheorem{conjecture}[theorem]{Conjecture}
\newtheorem{example}[theorem]{Example}

\newtheorem{definition}[theorem]{Definition}

\DeclareMathOperator{\dist}{dist}
\DeclareMathOperator{\Ker}{Ker}
\DeclareMathOperator{\ev}{ev}
\DeclareMathOperator{\ecc}{ecc}
\DeclareMathOperator{\spec}{sp}
\DeclareMathOperator{\supr}{sup}

\def\Q{\ns Q}

\def\Z{\ns Z}

\def\e{\mbox{\boldmath $e$}}

\def\v{\mbox{\boldmath $v$}}

\def\vecv{\mbox{\boldmath $v$}}

\def\vec0{\mbox{\boldmath $0$}}

\def\A{\mbox{\boldmath $A$}}
\def\B{\mbox{\boldmath $B$}}

\def\E{\mbox{\boldmath $E$}}
\def\D{\mbox{\boldmath $D$}}

\def\I{\mbox{\boldmath $I$}}

\def\L{\mbox{\boldmath $L$}}

\def\S{\mbox{\boldmath $S$}}

\def\U{\mbox{\boldmath $U$}}

\def\Z{\ns{Z}}
\def\I{\mbox{\boldmath $I$}}

\def\Q{\mbox{\boldmath $Q$}}

\def\1{\mbox{\boldmath $1$}}

\def\Re{\mathbb R}
\def\Z{\mathbb Z}

\begin{document}
	
	\title{
		On the spectra and spectral radii\\ of token graphs
		\thanks{The research of C. Dalf\'o and M. A. Fiol has been partially supported by AGAUR from the Catalan Government under project 2017SGR1087 and by MICINN from the Spanish Government under project PGC2018-095471-B-I00. 
			The research of M. A. Fiol was also supported by a grant from the  Universitat Polit\`ecnica de Catalunya with references AGRUPS-2022 and AGRUPS-2023.}
	}
	\author{M. A. Reyes$^a$, C. Dalf\'o$^a$, and  M. A. Fiol$^b$\\
		\\
		{\small $^a$Dept. de Matem\`atica, Universitat de Lleida, Igualada (Barcelona), Catalonia}\\
		{\small  \texttt{\{monicaandrea.reyes,cristina.dalfo\}@udl.cat}}\\
		{\small $^{b}$Dept. de Matem\`atiques, Universitat Polit\`ecnica de Catalunya, Barcelona, Catalonia} \\
		{\small Barcelona Graduate School of Mathematics} \\
		{\small  Institut de Matem\`atiques de la UPC-BarcelonaTech (IMTech)}\\
		{\small {\tt miguel.angel.fiol@upc.edu} }\\
	}

	\date{}
	\maketitle
	
	\begin{abstract}
		Let $G$ be a graph on $n$ vertices. The $k$-token graph (or symmetric $k$-th power) of $G$, denoted by $F_k(G)$  has as vertices the ${n\choose k}$  $k$-subsets of vertices from $G$, and two vertices are adjacent when their symmetric difference is a pair of adjacent vertices in $G$.
		In particular, $F_k(K_n)$ is the Johnson graph $J(n,k)$, which is a distance-regular graph used in coding theory.
		In this paper, we present some results concerning the (adjacency and Laplacian) spectrum of $F_k(G)$ in terms of the spectrum of $G$.
		For instance, when $G$ is walk-regular, an exact value for the spectral radius $\rho$ (or maximum eigenvalue) of $F_k(G)$ is obtained.
		When $G$ is distance-regular, other eigenvalues of its $2$-token graph are derived using the theory of equitable partitions.
		A generalization of Aldous' spectral gap conjecture (which is now a theorem) is proposed.
	\end{abstract}
	
	
	
	\section{Introduction}
	\label{sec:Intro}
	For a (simple and connected) graph $G=(V,E)$ with adjacency matrix $\A$,
	its local spectrum at vertex $u$ plays a role similar to the (standard adjacency) spectrum when the graph is `seen'
	from  vertex $u$. For instance, the local spectra of $G$, for every $u\in V$, were used by Fiol and Garriga \cite{fg97} to prove the so-called `spectral excess theorem', which gives a quasi spectral characterization of distance-regular graphs. In turn,
	this result was the crucial tool for the discovery, by van Dam and Koolen~\cite{vdk05}, of the first known family of non-vertex-transitive distance-regular graphs with unbounded diameter.
	Besides, Fiol, Garriga, and Yebra~\cite{fgy96}  used the local spectra to define the local predistance polynomials, which were used to characterize a general kind of local distance-regularity (intended for not necessarily regular graphs).
	
	One of the most important parameters in spectral graph theory is the \textit{index} or \textit{spectral radius} of a graph, which corresponds to the largest eigenvalue of its adjacency matrix.
	This parameter has special relevance in the study of many integer-valued graph invariants, such as the diameter, the radius, the domination number, the matching number, the clique number, the independence number, the chromatic number, or the sequence of vertex degrees. In turn, this leads to studying the structure of graphs having an extremal spectral radius and fixed values of some of such parameters. See Brualdi, Carmona, Van den Driessche, Kirkland, and Stevanovi\'c \cite[Cap. 3]{bcdks18}.
	
\begin{figure}[t]
		\begin{center}
			\includegraphics[width=6cm]{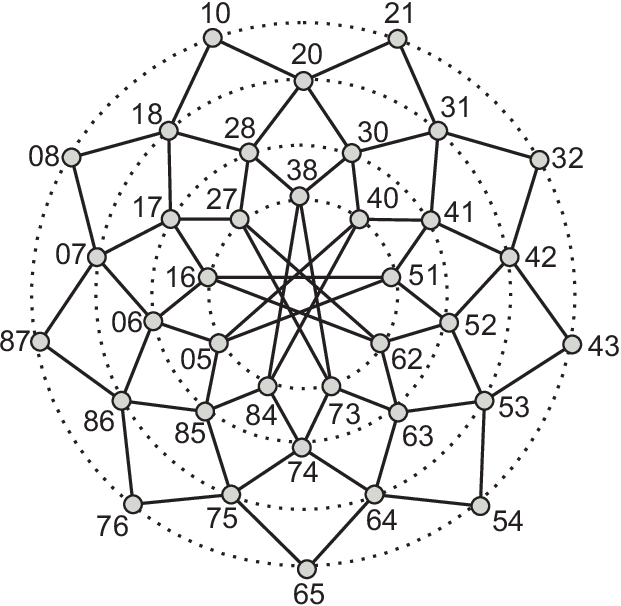}
			\caption{The $2$-token graph $F_2(C_9)$ of the cycle graph with vertex set $V(C_{9})=\Z_9$.
				The vertices inducing a circumference (in dashed line) of radius  $r_{\ell}$, with $\ell=1,2,3,4$ and $r_1>r_2>r_3>r_4$  are $ij$ with $\dist(i,j)=\ell$ in $C_9$.} 
			\label{fig1}
		\end{center}
	\end{figure}

	In this work, we use some information given by the local spectra to obtain new results about the spectral radius of an ample family of graphs, which are known as token graphs or symmetric $k$-th powers, defined as follows.
	For a given integer $k$, with $1\le k \le n$ (where $n$ is the order of $G$), the {\em $k$-token graph} $F_k(G)$ of $G$ is the graph whose vertex set $V (F_k(G))$ consists of the ${n \choose k}$
	$k$-subsets of vertices of $G$, and two vertices $A$ and $B$
	of $F_k(G)$ are adjacent whenever their symmetric difference $A \bigtriangleup B$ is a pair $\{a,b\}$ such that $a\in A$, $b\in B$, and $\{a,b\}\in E(G)$.
	In Figure \ref{fig1}, we show the 2-token graph of the cycle $C_9$ on 9 vertices.
	In particular,
	if $k=1$, then $F_1(G)\cong G$; and if $G$ is the complete graph $K_n$, then $F_k(K_n)\cong J(n,k)$, where $J(n,k)$ denotes the Johnson graph, see Fabila-Monroy,  Flores-Pe\~{n}aloza,  Huemer,  Hurtado,  Urrutia, and  Wood~\cite{ffhhuw12}.
	
	The name `token graph'
	also comes from the observation in
	\cite{ffhhuw12}, that vertices of $F_k(G)$ correspond to configurations
	of $k$ indistinguishable tokens placed at distinct vertices of $G$, where
	two configurations are adjacent whenever one configuration can be reached
	from the other by moving one token along an edge from its current position
	to an unoccupied vertex. Such graphs are also called {\em symmetric $k$-th power of a graph} in Audenaert, Godsil, Royle, and Rudolph \cite{agrr07}; and {\em $n$-tuple vertex graphs} in Alavi,  Lick, and Liu \cite{all02}. They have applications in physics; a connection between symmetric powers of graphs and the exchange of Hamiltonian operators in
	quantum mechanics is given in \cite{agrr07}.
	Our interest is in relation to the graph
	isomorphism problem. It is well known that there are cospectral non-isomorphic
	graphs, where often the spectrum of the adjacency matrix of a graph is used.
	For instance,
	Rudolph \cite{r02} showed that there are cospectral non-isomorphic graphs that can
	be distinguished by the adjacency spectra of their 2-token graphs, and he also gave an example for the Laplacian spectrum. Audenaert, Godsil, Royle, and Rudolph \cite{agrr07}
	also proved that 2-token graphs of strongly regular graphs with the same parameters
	are cospectral and also derived bounds on the (adjacency and Laplacian) eigenvalues of $F_2(G)$ for general graphs.
	For more information, see again 
	\cite{agrr07} or \cite{ffhhuw12}.
	
	What can be said about the spectrum of $F_k(G)$? 
	The three main results that we want to recall are the following.
	
	\begin{theorem}[Audenaert, Godsil, Royle, and  Rudolph \cite{agrr07}]
		All the strongly regular graphs with the same parameters have cospectral symmetric squares (or $2$-token graphs).
		\label{thm:Godsil}
	\end{theorem}
	
	\begin{theorem}[Dalf\'o, Duque, Fabila-Monroy,  Fiol, Huemer, Trujillo-Negrete,  Zaragoza Mart\'inez \cite{ddffhtz21}]
		\label{th:sp-token}
		For any graph $G$ on $n$ vertices, the Laplacian spectrum of its $h$-token graph is contained in the Laplacian spectrum of its $k$-token graph for every $1\le h\le k\le n/2$.
		\label{thm:DDFFHTZ}
	\end{theorem}
	
	\begin{theorem}[Lew \cite{l23}]
		\label{th:lew}
		Let $G$ have Laplacian eigenvalues $\lambda_1(=0)<\lambda_2\leq \cdots \leq\lambda_n$. Let $\lambda$ be an eigenvalue of $F_k(G)$ not in $F_{k-1}(G)$. Then,
		\begin{equation}
			\label{eq:lew}
			k(\lambda_2-k+1)\leq \lambda\leq k\lambda_n,
		\end{equation}
		and both bounds are tight.
	\end{theorem}
	
	In fact, the lower bound in \eqref{eq:lew} was first proved by Dalf\'o, Fiol, and Messegu\'e in  \cite[Ths. 3.5-3.6]{dfm22}.
	
	In this paper, we mainly derive new results about the spectral radius of token graphs, and it is organized as follows. The next section begins with some basic concepts, definitions, and results. More precisely, we recall some known results about the local spectra and derive the basic tools for computing the spectral radius.
	In Section \ref{sec:k-alpha}, we introduce the new concepts of $k$-algebraic connectivity and $k$-spectral radius. There, we study some of their properties and propose a generalization of Aldou's spectral gap conjecture, already a theorem (see Caputo, Liggett, and Richthammer \cite{clr10}). In Section \ref{sec:rho-token}, we give both lower and upper bounds for the spectral radius of a token graph, which are shown to be asymptotically tight.
	In the same section, we present some infinite families in which the exact values of the spectral radius are obtained.
	Finally, in the last section, we deal with the case of distance-regular and strongly regular graphs, where two results are presented in the form of Audenaert, Godsil, Royle, and Rudolph's result \cite{agrr07}, and Lew's result \cite{l23}.
	
	\section{Preliminaries}
	\label{sec:Prelim}
	
	\subsection{Graphs and their spectra}
	
	Let $G$ be a (simple and connected) graph with vertex set $V(G)=\{1,2,\ldots,n\}$ and edge set $E(G)$.
	Let $G$ have adjacency matrix $\A$, and spectrum
	$$
	\spec G\equiv \spec \A=\{\theta_0^{m_0},\theta_1^{m_1},\ldots,\theta_d^{m_d}\},
	$$
	where $\theta_0>\theta_1>\cdots >\theta_d$.
	Thus, by the Perron-Frobenius theorem, $G$ has spectral radius $\rho(G)=\theta_0$.
	
	Let $\L=\D-\A$ be the Laplacian matrix of $G$, with eigenvalues
	$$
	\lambda_1(=0)<\lambda_2\le \cdots \le\lambda_n.
	$$
	Recall that $\lambda_2$ is the algebraic connectivity, and  $\D$ is the diagonal matrix whose diagonal entries are the vertex degrees of $G$.
	
	\subsection{The local spectra of a graph}
	\label{subsec:local-sp}
	
	Let $G$ have different eigenvalues $\theta_0>\cdots>\theta_d$, with respective multiplicities $m_0,\ldots,m_d$.
	If $\U_i$ is the $n\times m_i$ matrix whose columns are the orthonormal eigenvectors of $\theta_i$,
	the matrix $\E_i=\U_i\U_i^{\top}$, for $i=0,1,\dots,d$, is the {\it $($principal\/$)$ idempotent} of
	$\A$ and represents the orthogonal projection of $\Re^n$
	onto the eigenspace $\Ker (\A-\theta_i \I)$.
	The {\em
		$(u$-$)$local multiplicities} of the eigenvalue $\theta_i$ are
	defined as
	$$
	m_u(\theta_i) = \|\E_i\e_u\|^2
	= \langle\E_i\e_u,\e_u\rangle =
	(\E_i)_{uu}\qquad \mbox{for } u\in V \mbox{ and  } i =0,1,\dots,d,
	$$
	where $\e_u$ is the $n$-dimensional vector with a $1$ in the $u$-th entry and zeros elsewhere. In particular, $m_u(\theta_0)=v_{u}^2>0$, where $\v$ is the corresponding normalized Perron eigenvector, that is, the eigenvector that can be chosen to have strictly positive components.
	Although the local multiplicities are, of course, not necessarily integers, they have nice properties when the graph is studied from a vertex, so justifying their name. Thus,
	they satisfy $\sum_{i=0}^d m_u(\theta_i) = 1$ and $\sum_{u\in
		V} m_u(\theta_i) =m_i$, for $i=0,1,\dots,d$. The number $a_{uu}^{(\ell)}$ of closed walks of
	length $\ell$ rooted at vertex $u$ can be computed as
	\begin{equation}
		\label{closed-walks}
		a_{uu}^{(\ell)}=\sum_{i=0}^d m_{u}(\theta_i) \theta_i^{\ell}
	\end{equation}
	(see Fiol, Garriga, and Yebra \cite[Corollary 2.2]{fgy96}).
	By picking up the eigenvalues $\theta_i$ with non-null  local multiplicities,
	$\mu_0(=\theta_0)>\mu_1>\cdots>\mu_{d_u}$, 
	we define the {\it
		($u$-)local spectrum} of $G$ as
	$$
	\spec_u G:=
	\{\mu_0^{m_{u}(\mu_0)},\mu_1^{m_{u}(\mu_1)},\ldots,\mu_{d_u}^{m_{u}(\mu_{d_u})}\},
	$$
	with {\it ($u$-)local mesh}, or set of distinct eigenvalues,
	$\ev_u G:=\{\mu_0>\mu_1>\cdots>\mu_{d_u}\}$.
	The eccentricity of a vertex $u$ satisfies an upper bound similar to that satisfied by the diameter of $G$ in terms of its distinct eigenvalues. More precisely,
	\begin{equation}
		\label{bound-ecc}
		\ecc (u) \leq d_u=|\ev_u G|-1.
	\end{equation}
	In coding theory, $d_u$ corresponds to the
	so-called `dual degree' of the trivial code $\{u\}$.
	For more information, see Fiol, Garriga, and Yebra \cite{fgy96}.
	
	We use the following lemma to prove the results of Section \ref{sec:rho-token}. Notice that this is just a reformulation of the power method in terms of the number of walks given by \eqref{closed-walks}.
\begin{lemma}
	\label{basic-lemma}
	Let $G$ be a finite graph with different eigenvalues $\theta_0>\cdots > \theta_d$. Let $a^{(\ell)}_u$ be the number of $\ell$-walks starting from (any fixed) vertex $u$, and let $a^{(\ell)}_{uu}$ be the number of closed $\ell$-walks rooted at $u$. Then,
	$$
	\rho(G)=\lim_{\ell \rightarrow \infty} \sqrt[\ell]{a^{(\ell)}_u}=\lim_{\ell \rightarrow \infty} \supr\sqrt[\ell]{a^{(\ell)}_{uu}},
	$$
	where $`\supr$' denotes the supremum.
\end{lemma}
\begin{proof}
	We only prove that $\rho(G)$ equals the second limit (the other equality proved similarly). Using \eqref{closed-walks},
	\begin{align*}
		a_{uu}^{(\ell)} & =\sum_{i=0}^d m_{u}(\theta_i) \theta_i^{\ell}
		= \theta_0^{\ell}\left(m_{u}(\theta_0) + m_{u}(\theta_1)\left(\frac{\theta_1}{\theta_0}\right)^{\ell} +\cdots
		+ m_{u}(\theta_d)\left(\frac{\theta_d}{\theta_0}\right)^{\ell}\right)\\
		& \stackrel{\ell \rightarrow \infty}{\longrightarrow}
		\left\{
		\begin{array}{cl}
			2\theta_0^{\ell} m_u(\theta_0) & \mbox{if $G$ is bipartite and $\ell$ is even},  \\
			0 & \mbox{if $G$ is bipartite and $\ell$ is odd,}\\ [.1cm]
			\theta_0^{\ell}m_u(\theta_0) & \mbox{otherwise,}
		\end{array}
		\right.
	\end{align*}
	where we used that $m_u(\theta_0)>0$ and $\theta_0>|\theta_i|$ for every $i\neq 0$, except when $G$ is bipartite, in which case $\theta_d=-\theta_0$ and $m_u(\theta_d)=m_u(\theta_0)$.
	Thus, the result holds by taking $\ell$-th roots.
\end{proof}

\subsection{Regular partitions and their spectra}
\label{subsec:reg-part}

Dealing with a regular partition is a useful method in graph theory, as it allows us to obtain some information about a graph considering a smaller version of it (the so-called `quotient graph'). Besides, if we consider a graph $G$ and a group of automorphisms $\Gamma$, then a partition of the vertices of the graph into orbits by $\Gamma$ is a regular partition.

Let $G=(V,E)$ be a graph with vertex set $V=V(G)$, adjacency matrix $\A$, and Laplacian matrix $\L$. A partition $\pi$ of its vertex set $V$ into $r$ cells $C_1,C_2,\ldots,C_r$ is called {\em regular} (or {\em equitable})
whenever, for any $i,j=1,\ldots,r$, the {\em intersection numbers}
$b_{ij}(u)=|G(u)\cap C_j|$ (where $u\in C_i$, and $G(u)$ is the set of vertices that are neighbors of the vertex $u$) do not depend on the vertex $u$ but only on the cells $C_i$ and $C_j$. In this case, such numbers are simply written as $b_{ij}$, and the $r\times r$ matrices $\Q_A=\A(G/\pi)$ and $\Q_L=\L(G/\pi)$ with entries $(\Q_A)_{ij}=b_{ij}$ and
\begin{equation}
	\textstyle
	(\Q_L)_{ij}=\left\{
	\begin{array}{cl}
		-b_{ij} & \mbox{if } i\neq j,\\
		b_{ii}-\displaystyle \sum_{j=1}^r b_{ij} & \mbox{if } i=j ,
	\end{array}
	\right.
	\label{Q_L}
\end{equation}
are, respectively, referred to as the \textit{quotient matrix} and \textit{quotient Laplacian matrix} of $G$ with respect to $\pi$.
In turn, these matrices correspond to the \textit{quotient (weighted) directed graph} $G/\pi$,
whose vertices representing the $r$ cells, and there is an arc with weight $b_{ij}$ from vertex $C_i$ to vertex $C_j$ if and only if $b_{ij}\neq 0$.
Of course, if $b_{ii}>0$, for some $i=1,\ldots,r$, the quotient graph (or digraph) $G/\pi$ has loops.
Given a partition $\pi$ of $V$ with $r$ cells, let $\S$ be the {\em characteristic matrix} of $\pi$, that is, the $n\times r$ times matrix whose columns are the characteristic vectors of the cells of $\pi$. Then, $\pi$ is a regular partition if and only if $\A\S=\S\Q_A$ or $\L\S=\S\Q_L$. Moreover,
$$
\Q_A=(\S^{\top}\S)^{-1}\S^{\top}\A\S,\quad\mbox{and}\quad \Q_L=(\S^{\top}\S)^{-1}\S^{\top}\L\S.
$$
Thus, there is a strong analogy with similar results satisfied by the Laplacian matrices of the $h$-token graph and $k$-token graph of $G$ for $h\le k$. More precisely, we have the following statements (Those numbered with 1 can be found in Godsil \cite{g93,g95}. The statements in 2 are derived similarly to those in number 1 when we use the quotient Laplacian matrix, defined in  \eqref{Q_L}, instead of the standard quotient matrix.  Finally, the statements in  3 follow from the results given by Dalf\'o, Duque, Fabila-Monroy, Fiol, Huemer, Trujillo-Negrete in \cite{ddffhtz21}.):
\begin{enumerate}
	\item
	If $\pi$ is a regular partition with characteristic matrix $\S$ and quotient matrix $\Q_A=\A(G/\pi)$, then
	\begin{enumerate}
		\item
		$\A\S=\S\Q_A$.
		\item
		$\Q_A=(\S^{\top}\S)^{-1}\S^{\top}\A\S$.
		\item
		The column space (and its orthogonal complement) of $\S$ is $\A$-invariant.
		\item
		The characteristic polynomial of $\Q_A$ divides the characteristic polynomial of $\A$. Thus, $\spec \Q_A \subseteq \spec \A$.
		\item
		If $\v$ is an eigenvector of $\Q_A$ with eigenvalue $\lambda$, 
		then $\S\v$ is an eigenvector of $\A$ with the same eigenvalue. 
		(The eigenvector $\v$ of $\Q_A$ {\em `lifts'} to an eigenvector of $\A$.)
		\item
		If $\v$ is an eigenvector of $\A$ with eigenvalue $\lambda$ and $\S^{\top}\v\neq\vec0$,
		then $\S^{\top}\v$ is an eigenvector of $\Q_A$ with the same eigenvalue.
	\end{enumerate}
	\item
	If  $\pi$ is a regular partition with characteristic matrix $\S$ and quotient Laplacian matrix $\Q_L=\L(G/\pi)$, then
	\begin{enumerate}
		\item
		$\L\S=\S\Q_L$.
		\item
		$\Q_L=(\S^{\top}\S)^{-1}\S^{\top}\L\S$.
		\item
		The column space (and its orthogonal complement) of $\S$ is $\L$-invariant.
		\item
		The characteristic polynomial of $\Q_L$ divides the characteristic polynomial of $\L$.
		Thus, $\spec \Q_L \subseteq \spec \L$.
		\item
		If $\v$ is an eigenvector of $\Q_L$ with eigenvalue $\lambda$,
		then $\S\v$ is an eigenvector of $\L$ with the same eigenvalue.
		\item
		If $\v$ is an eigenvector of $\L$ with eigenvalue $\lambda$ and $\S^{\top}\v\neq\vec0$,
		then $\S^{\top}\v$ is an eigenvector of $\Q_L$ with the same eigenvalue.
	\end{enumerate}

	\item
	Let $\L_h=\L(F_h(G))$ and $\L_k=\L(F_k(G))$ be, respectively,  the Laplacian matrices of the $h$-token and $k$-token graphs of $G$, for $h\le k$, and  let $\S_b$ be the  $(k,h)$-binomial matrix.
	This is a ${n \choose k}\times{n \choose h}$ matrix whose rows are indexed by the $k$-subsets of $A\subset [n]$, and its columns are indexed by the $h$-subsets of $X\subset [n]$, with entries
	$$
	(\S_b)_{AX}=
	\left\lbrace
	\begin{array}{ll}
		1 & \mbox{if } X\subset A,\\
		0 & \mbox{otherwise.}
	\end{array}
	\right.
	$$
	Then,
	\begin{enumerate}
		\item
		$\L_k\S_b=\S_b\L_h$.
		\item
		$\L_h=(\S_b^{\top}\S_b)^{-1}\S_b^{\top}\L_k\S_b$.
		\item
		The column space (and its orthogonal complement) of $\S_b$ is $\L_k$-invariant.
		\item
		The characteristic polynomial of $\L_h$ divides the characteristic polynomial of $\L_k$.
		Thus, $\spec \L_h \subseteq \spec \L_k$.
		\item
		If $\v$ is an eigenvector of $\L_h$ with eigenvalue $\lambda$,
		then $\S_b\v$ is an eigenvector of $\L_k$ with eigenvalue $\lambda$.
		\item
		If $\v$ is an eigenvector of $\L_k$ with eigenvalue $\lambda$ and $\S_b^{\top}\v\neq\vec0$,
		then $\S_b^{\top}\v$ is an eigenvector of $\L_h$ with the same eigenvalue.
	\end{enumerate}
\end{enumerate}

\subsection{Walk-regular graphs}
\label{subsec:walk-reg}

Let $a_u^{({\ell})}$ denote the number of closed walks of
length ${\ell}$ rooted at vertex $u$, that is,
$a_u^{({\ell})}=a_{uu}^{({\ell})}$. If these numbers only
depend on ${\ell}$, for each $\ell \ge 0$, then $G$ is called
{\em walk-regular}, a concept introduced by Godsil and
McKay in \cite{gmk80}.

Notice that, as
$a_u^{(2)}=\delta_u$, the degree of vertex $u$, a walk-regular graph is necessarily regular.

Moreover,
a graph $G$ is called {\it spectrally regular} when all vertices have the same local spectrum:
$\spec_u G=\spec_v G$ 
for any $u,v\in V$. The following result (in Delorme and Tillich \cite{dt97}, Fiol and Garriga \cite{fg98}, and also Godsil and McKay \cite{gmk80}) provide some characterizations of such graphs. 

\begin{lemma}[\cite{dt97},\cite{fg98},\cite{gmk80}]
	Let $G=(V,E)$ be a graph. The following statements are equivalent.
	\begin{itemize}
		\item[$(i)$]
		$G$ is walk-regular.
		\item[$(ii)$]
		$G$ is spectrally regular.
		\item[$(iii)$]
		The spectra of the vertex-deleted subgraphs are all equal: $\spec\,
		(G\setminus u)=\spec\, (G\setminus v)$ for any $u,v\in V$.
	\end{itemize}
\end{lemma}

\section{The $k$-algebraic connectivity and $k$-spectral radius}
\label{sec:k-alpha}

In this section, we always consider the Laplacian spectrum. Let $G$ be a graph on $n$ vertices, and $F_k(G)$ its $k$-token graph for $k\in \{0,1,\ldots,n\}$. Note that $F_k(G)\cong F_{n-k}(G)$ where, by convenience, $F_0(G)\cong F_n(G)=K_1$ (a singleton). Moreover, $F_1(G)\cong G$. From 
Dalf\'o, Duque, Fabila-Monroy, Fiol, Huemer, Trujillo-Negrete, and Zaragoza Mart\'{\i}nez \cite{ddffhtz21}, it is known that the Laplacian spectra of the token graphs of $G$ satisfy
\begin{equation}
	\label{sp-inclusion}
	\{0\}=\spec F_0(G)\subset \spec F_1(G)\subset \spec F_2(G)\subset \cdots \subset \spec F_{\lfloor n/2\rfloor}(G).
\end{equation}
Let denote $\alpha(G)$ and $\rho(G)$ the algebraic connectivity (see Fiedler \cite{fi73}) and the spectral radius of a graph $G$, respectively. Then, from \eqref{sp-inclusion}, we have
\begin{align}
	\alpha(G) &\ge \alpha(F_2(G))\ge \cdots \ge  \alpha(F_{\lfloor n/2\rfloor}(G)), \label{ineq-alpha's}\\
	\rho(G) &\le \rho(F_2(G))\le \cdots \le \rho(F_{\lfloor n/2\rfloor}(G)). \label{ineq-rho's}
\end{align}
The concepts of algebraic connectivity  and spectral radius, together with \eqref{sp-inclusion}--\eqref{ineq-rho's}, suggest the following definitions.
\begin{definition}
	Given a graph $G$ on $n$ vertices and an integer $k$ such that $1\le k\le \lfloor n/2\rfloor$, the {\em $k$-algebraic connectivity} $\alpha_k=\alpha_k(G)$ and the {\em $k$-spectral radius} $\rho_k=\rho_k(G)$  of $G$ are, respectively, the minimum and maximum eigenvalues of the multiset $\spec F_k(G)\setminus \spec F_{k-1}(G)$.
\end{definition}

Notice that, since ${n\choose k}>{n\choose k-1}$ for $1\le k\le \lfloor n/2\rfloor$, the parameters $\alpha_k$ and $\rho_k$ always exist.

For instance, with $G=P_6$, the path on $6$ vertices, we have (approximately) 
$$
\alpha_1(P_6)=0.2679,\quad \alpha_2(P_6)=0.5727,\quad
\alpha_3(P_6)=0.9302, 
$$
and
$$
\rho_1(P_6)=3.7320,\quad  \rho_2(P_6)=6.504,\quad\rho_3(P_6)=7.4871,
$$
whereas for $G=C_7$, the cycle on $7$ vertices, we get
$$
\alpha_1(C_7)=0.7530,\quad  \alpha_2(C_7)=1.1633,\quad
\alpha_3(C_7)=1.2696. 
$$

Moreover, from these definitions, the following facts hold.
\begin{enumerate}
	\item[$(i)$] 
	$\rho_k(G)\ge \alpha_k(G)\ge 0$.
	\item[$(ii)$]
	$\alpha_1(G)=\alpha(G)$ (the standard algebraic connectivity of $G$) and $\rho_1(G)=\rho(G)$ (the standard spectral radius of $G$).
	\item[$(iii)$]
	Since $F_k(K_n)\cong J(n,k)$ (the Johnson graph), we have
	$$
	\alpha_k(K_n)=\rho_k(K_n)=k(n+1-k), \qquad k=1,\ldots,\lfloor n/2\rfloor.
	$$
	In particular, $\alpha_1(K_n)=\rho_1(K_n)=n$, $\alpha_2(K_n)=\rho_2(K_n)=2(n-1)$, and so on. 
\end{enumerate}
The equalities in $(iii)$ come from the fact that the Johnson graph $J(n,k)$ has different Laplacian eigenvalues $\lambda_j=j(n+1-k)$, with multiplicities $m_j={n\choose j}-{n\choose j-1}$ for $j=0,1,\ldots,k$.

From what is known about token graphs, we can suggest
some conjectures and state some results, as follows.
\begin{conjecture}
	\label{conjec-creix-alfa}
	For any graph $G$,
	$$
	\alpha_1(G)\le \alpha_2(G)\le \cdots \le \alpha_{\lfloor n/2\rfloor}(G).
	$$
\end{conjecture}

Because of \eqref{ineq-alpha's}, if Conjecture \ref{conjec-creix-alfa} holds, also does the conjecture proposed in \cite{ddffhtz21}, that is,  $\alpha(F_k(G))=\alpha(G)$ for any $k\le n/2$.
In fact, the last equality follows from the proof of Aldous' spectral gap conjecture given by Caputo, Ligget, and Richthammer in \cite{clr10}. By this result, what we can state is that
$
\min\{\alpha_2,\ldots, \alpha_{\lfloor n/2\rfloor}\}\ge \alpha_1.
$

\begin{conjecture}
	For any graph $G$,
	$$
	\rho_1(G)\le \rho_2(G)\le \cdots \le \rho_{\lfloor n/2\rfloor}(G).
	$$
\end{conjecture}
Notice that, from \eqref{ineq-rho's}, if this conjecture holds, then 
$\rho_k(G)=\rho(F_k(G))$ for any $k\le n/2$.
\begin{lemma}
	\label{lem:alpha+rho}
	For any graph $G$ and its complementary graph $\overline{G}$, the $k$-algebraic connectivity and $k$-spectral radius of $\overline{G}$ satisfy
	$$
	\alpha_k(G)+\rho_k(\overline{G})=k(n-k+1).
	$$
\end{lemma}
\begin{proof}
	It was proved in \cite{ddffhtz21} that each eigenvalue of $J(n,k)$ is the sum of one eigenvalue of $F_k(G)$ and one eigenvalue of $F_k(\overline{G})$. Then, from 
	\eqref{sp-inclusion}, the eigenvalues of $\spec F_k(G)\setminus \spec F_{k-1}(G)$
	and $\spec F_k(\overline{G})\setminus \spec F_{k-1}(\overline{G})$ must be paired in such a way that their sums are all equal to the `new' eigenvalue of $\spec F_k(J(n,k))\setminus \spec F_{k-1}(J(n,k))$, which is $k(n-k+1)$. In particular, this happens with the minimum eigenvalue of $\spec F_k(G)\setminus \spec F_{k-1}(G)$, $\alpha_k(G)$ and the maximum eigenvalue of $\spec F_k(\overline{G})\setminus \spec F_{k-1}((\overline{G})$, $\rho_k(\overline{G})$.
\end{proof}

Moreover, $\alpha(G)=n-\rho(\overline{G})$, as it is well known.
\begin{corollary}
	For any graph $G$ on $n$ vertices,
	\begin{align*}
		&\alpha_k(G)\le k(n-k+1),\qquad k=1,\ldots, \lfloor n/2\rfloor.\\
		&\rho_k(\overline{G})\le k(n-k+1),\qquad k=1,\ldots, \lfloor n/2\rfloor.\\
	\end{align*}
\end{corollary}

From Lemma \ref{lem:alpha+rho} and Proposition \ref{propo:sp-radius}, we get the following result, which will be proved in Section \ref{sec:drg's}.
\begin{corollary}
	Let $G$ be a bipartite distance-regular graph.
	Let  $\L(F_2/\pi)$ be  the quotient matrix in \eqref{L-quotient} with spectral radius  $\rho_L(F_2/\pi)$.
	Then, 
	$$
	\alpha_2(\overline{G})={n\choose 2}-\rho_L(F_2/\pi).
	$$
\end{corollary}


\section{The spectral radius of token graphs}
\label{sec:rho-token}

In contrast with the previous section, in this section, we always consider the spectral radius of the adjacency matrix of a (connected) graph. Consider a graph $G$ with spectral radius $\rho(G)$ and vertex-connectivity $\kappa$ (that is, the minimum number of vertices whose suppression either disconnects the graph or results in a singleton).
By taking the spectral radii of its $U$-deleted subgraphs, with $U\subset V$ and $|U|=k<\kappa$, we define the following two parameters:
\begin{align*}
	\rho^k_M(G) &=\max\{\rho(G\setminus U): U\subset V, \,|U|=k\},\\
	\rho^k_m(G) &=\min\{\rho(G\setminus U):  U\subset V, \,|U|=k\}.
\end{align*}
Notice that, if $G$ is walk-regular, then $\rho^1_M(G)=\rho^1_m(G) =\rho(G\setminus u)$ for every vertex $u$.
If $G$ is distance-regular with degree $\delta$, it is known that it has vertex-connectivity $\kappa(G)=\delta$ (see Brouwer and  Koolen \cite{bk09}). Moreover,  Dalf\'{o}, Van Dam, and Fiol \cite{dvdf11} showed that $\spec (G\setminus U)$ only depends on the distances in $G$ between the vertices of $U$.
Thus,  for every $k\le \delta-1$, the computation of $\rho^k_M(G)$ and $\rho^k_m(G)$ can be drastically reduced by considering only the subsets $U$ with different `distance-pattern' between vertices. For instance, if $G$ has diameter $D$,
\begin{align*}
	\rho^2_M(G) &=\max_{1\le \ell\le D}\{\rho(G\setminus \{u,v\}): \dist_G(u,v)=\ell\},\\
	\rho^2_m(G) &=\min_{1\le \ell\le D}\{\rho(G\setminus \{u,v\}): \dist_G(u,v)=\ell\}.
\end{align*}

In general, by using interlacing (see Haemers \cite{h95} or Fiol \cite{f99}), we have the following result.
\begin{lemma}
	Let $G$ be a graph with $n$ vertices, vertex-connectivity $\kappa$, and eigenvalues $\lambda_1\ge \lambda_2\ge \cdots \ge \lambda_n$.
	Then, for every $k=1,\ldots,\kappa-1$, 
	\begin{align}
		\lambda_{k+1} & \le \rho^k_M(G)\le \lambda_1,\\
		\lambda_{n} & \le \rho^k_m(G)\le \lambda_{n-k}.
	\end{align}
\end{lemma}

From the above results, Lemma \ref{basic-lemma}, and the bounds for the spectral radius of graph perturbations obtained in Dalf\'o, Garriga, and Fiol \cite{dfg11} and Nikiforov \cite{nk07}, we obtain the following  main result.

\begin{theorem}
	Let $G$ be a graph with spectral radius $\rho(G)$ and vertex-connectivity $\kappa>1$. Given an integer $k$, with $1\le k<\kappa$, let $\rho^k_M(G)$ and  $\rho^k_m(G)$ be the maximum and minimum of the spectral radii of the $U$-deleted subgraphs of $G$, where $|U|=k$. 
	\begin{itemize}
		\item[$(i)$]
		The spectral radius of the $k$-token graph $F_k(G)$ satisfies
		\begin{equation}
			\label{basin-ineq}
			k\rho^{k-1}_m(G)\le \rho(F_k(G))\le k\rho^{k-1}_M(G).
		\end{equation}
		\item[$(ii)$]
		If $G$ is a graph of order $n$ and diameter $D$, the spectral radius of the $k$-token graph $F_k(G)$ satisfies
		\begin{equation}
			\label{basin-ineq2}
			\rho(F_k(G))< k\left(\rho(G)-\frac{1}{n\rho(G)^{2D}}\right).
		\end{equation}
		\item[$(iii)$]
		If $G$ is walk-regular and $k=2$ (that is, $F_2(G)$ is the 2-token graph of $G$), then 
		\begin{equation}
			\rho(F_2(G))= 2\rho^1_m(G)= 2\rho^1_M(G).
		\end{equation}
	\end{itemize}
\end{theorem}
\begin{proof}
	$(i)$ To prove the upper bound in \eqref{basin-ineq}
	the key idea is the following: Given some integer $\ell$ large enough,
	every $\ell$-walk $W$ in $F_k(G)$ from a given vertex $A$ can be seen as $\ell_i$-walks in $G$ for $i=1,\ldots, k$ such that $\sum_{i=1}^k \ell_i=\ell$. Each step of $W$ corresponds to one step given by a token, say `1', whereas the remaining tokens are `still'. Thus, the move of token `1' is done in $G\setminus U$ for some vertex subset $U$, with $|U|=k-1$, and the same holds for every token. Thus, for $\ell$ great enough, 
	$W$ induces a number of walks of order ${\ell\choose \ell_1,\ell_2,\ldots, \ell_k}$ since, in each step, one of the $k$ tokens can be moved. For large $\ell$, the number of `conflicts' (that is, two 
	moves of tokens to the same vertex of $G$) is negligible with respect to the length of the $\ell$-walk.
	Then, the total number of $\ell$-walks $w_A^{(\ell)}$ in $F_k(G)$ starting from $A$ is of order 
	$$
	w_A^{(\ell)}\sim \sum_{\ell_1+\cdots+\ell_k=\ell}{\ell\choose \ell_1,\ell_2,\ldots, \ell_k}w_u^{(\ell)}=k^{\ell}w_u^{(\ell)},
	$$
	(in which we applied the Multinomial Theorem to obtain the equality) and  where, for some given vertex $u$ in $G$, $w_u^{(\ell)}$ is the number of $\ell$-walks starting from $u$ in $G\setminus U$.
	Consequently
	$$
	\rho(F_k(G))=\lim_{\ell\rightarrow\infty}\sqrt[\ell]{w_A^{(\ell)}}=k \lim_{\ell\rightarrow\infty}\sqrt[\ell]{w_u^{(\ell)}}\le k\rho_M^{k-1}(G).
	$$
	The lower bound comes from $\lim_{\ell\rightarrow\infty}\sqrt[\ell]{w_u^{(\ell)}}\ge k\rho_M^{k-1}(G)$.
	
	$(ii)$ This follows from $(i)$ and a result of Nikiforov (see \cite[Theorem 1]{nk07b}) stating that if $G$ is a graph as in the conditions of the theorem, and $H(=G\setminus U)$ is a proper subgraph of $G$, then $\rho(H)<\rho(G)-\frac{1}{n\rho(G)^{2D}}$.\\
	Finally, $(iii)$ holds by using \eqref{basin-ineq} since, as commented above, in a walk-regular graph $G$, we have $\rho^1_M(G)=\rho^1_m(G) =\rho(G\setminus u)$ for every vertex $u$.
\end{proof}

Notice that, if $G$ is connected $(D<\infty$), \eqref{basin-ineq2} improves the upper bound $k \rho(G)$ of Lew \cite{l23} in Theorem \ref{th:lew}.
Moreover,  as commented in Nikiforov \cite{nk07b}, for large values of $\rho(G)$ and $D$, the right hand of \eqref{basin-ineq2} yields the correct order of magnitude of $\rho(H)$. Thus, we can say that, asymptotically, the spectral radius of $F_k(G)$ is $k$ times the spectral radius of $G$. 
Moreover, in the case when $G$ is regular, $(ii)$  can be rewritten as
\begin{equation}
	\rho(F_{k}(G))<  k\left(\rho(G)-\frac{1}{n(D+1)}\right),
\end{equation}
(see again \cite[Th. 4]{nk07b}).

Since the different eigenvalues of the path $P_n$ on $n$ vertices are $\theta_i=2\cos\left(\frac{i\pi}{n+1} \right)$ for $i=1,\ldots,n$,
and the spectral radius of the complete bipartite graph is $\rho(K_{m,n})=\sqrt{mn}$, we get the following results.

\begin{corollary}
Let $P_n$ and $C_n$ be, respectively, the path and cycle graph on $n$ vertices. 
Let $P_{\infty}$ and $C_{\infty}$ be, respectively, the infinite path and cycle graphs. Then, their spectral radii satisfy the following:
\begin{itemize}
	\item[$(i)$]
	$\rho(F_2(P_n))\le 4\cos(\pi/n)$ and $\rho(F_2(P_{\infty}))= 4$.
	\item[$(ii)$]
	$\rho(F_2(C_n))= 4\cos(\pi/n)$ and $\rho(F_2(C_{\infty}))= 4$.
	\item[$(iii)$]
	$\rho(F_2(K_{n,n}))= 2\sqrt{n(n-1)}$.
\end{itemize}
\end{corollary}

\begin{table}[!ht]
\begin{center}
	\begin{tabular}{|c|ccccccc| }
		\hline
		$n$ & $3$ & $4$ & $\cdots$  & $8$ & $9$ & $10$ & $11$ \\
		\hline\hline
		$\rho(P_{n-1})$ & 1 & 1.41421 & $\cdots$  & 1.84776 & 1.87938 & 1.92113& 1.91898 \\
		\hline
		$\rho(F_2(C_n))$  & 2 & 2.82842 & $\cdots$ & 3.69552  & 3.75877 & 3.84226 & 3.83796 \\
		\hline
	\end{tabular}
\end{center}
\caption{Spectral radii of the 2-token graphs of the cycles $C_n$ with respect to spectral radii of the paths graphs $P_{n-1}$.}
\label{table4}
\end{table}

Let us show a pair of examples:
\begin{itemize}
\item
$F_2(C_3)=C_3=K_3$ has spectrum $\{-1^{[2]},2\}$, whereas $P_2$ has $\{-1,1\}$.
\item
$F_2(C_4)=K_{2,4}$ has spectrum $\{-2\sqrt{2},0^{[6]},2\sqrt{2}\}$,  whereas $P_3$ has $\{-\sqrt{2},0,\sqrt{2}\}$.
\end{itemize}


\section{The case of distance-regular graphs}
\label{sec:drg's}

In this section, we adopt the terminology of 
Brouwer, Cohen, and Neumaier \cite{bcn89} for distance-regular graphs. Furthermore, since we examine both the adjacency and Laplacian spectra, we  denote their respective spectral radii as $\rho_A$ and $\rho_L$. In the following result, consider that $G$ is a distance-regular graph with degree $\delta=b_0$, diameter $d$, intersection array
\begin{equation}
\label{intersec-array}
\iota(G)=\{
b_0, b_1, \ldots, b_{d-1};
c_1, c_2, \ldots,  c_d\}.
\end{equation}
or intersection matrix
\begin{align}
\B & =\left(
\begin{array}{ccccc}
0   & b_0 &        &        &      \\
c_1 & a_1 & b_1    &        &      \\
& c_2 & a_2    & \ddots &      \\
&     & \ddots & \ddots & b_{d-1} \\
&     &        & c_d    &  a_d \\
\end{array}
\right),\label{A-quotient}
\end{align}
where $a_i=\delta-b_i-c_i$, for $i=1,\ldots,d$.

\begin{lemma}
\label{regp-drg}
Let $F_2(G)$ be the 2-token graph of a distance-regular graph $G$ 
with degree $\delta=b_0$, diameter $d$, and intersection array $\iota(G)$ as in \eqref{intersec-array}. 
Then, $F_2=F_2(G)$ has a regular partition $\pi$ with quotient matrix and  quotient Laplacian matrix
Then, $F_2=F_2(G)$ has a regular partition $\pi$ with quotient matrix and  quotient Laplacian matrix
\begin{align}
\A(F_2/\pi) & =2\left(
\begin{array}{ccccc}
	a_1 & b_1 &        &        &      \\
	c_2 & a_2 & b_2    &        &      \\
	& c_3 & a_3    & \ddots &      \\
	&     & \ddots & \ddots &   b_{d-1}\\
	&     &        & c_d    &  a_d \\
\end{array}
\right),\label{A-quotient2}\\
\L(F_2/\pi) & =2\left(
\begin{array}{ccccc}
	b_1  & -b_1    &         &         &      \\
	-c_2 & c_2+b_2 & -b_2    &         &      \\
	& -c_3    & c_3+b_3 & \ddots  &      \\
	&         & \ddots  & \ddots  &  -b_{d-1}  \\
	&         &         & -c_d    &  c_d \\
\end{array}
\right), \label{L-quotient}
\end{align}
where $c_i+a_i+b_i=\delta$, for $i=0,1,\ldots,d$.
\end{lemma}
\begin{proof}
Let us consider the partition $\pi$ with classes $C_1,C_2,\ldots,C_d$, where $C_i$ is constituted by the vertices $\{u,v\}$ of $F_2(G)$, with $u,v\in V(G)$, such that $\dist_G(u,v)=i\in \{1,\ldots,d\}$. Note that $\A(F_2/\pi)$ is a $d\times d$ matrix since $\pi$ has $d$ classes. Then, since $G(v)\cap G_{i-1}(u)=c_{i-1}$, $G(v)\cap G_{i}(u)=a_i$ and $G(v)\cap G_{i+1}(u)=b_{i+1}$ (and the same equalities hold when we interchange $u$ and $v$), the vertex  $\{u,v\}\in C_i$ is adjacent to $2c_{i-1}$, $2a_i$, and 2$b_{i+1}$ vertices in $C_{i-1}$, $C_i$, and $C_{i+1}$, respectively. Thus, since this holds for every vertex in $C_i$, the partition $\pi$ is regular with the quotient matrix in \eqref{A-quotient2}. From this and \eqref{Q_L}, we obtain the Laplacian quotient matrix of 
\eqref{L-quotient}.
\end{proof}

The following result shows that the quotient matrices of a regular partition can be used to find the spectral radius or Laplacian spectral radius of the $2$-token graph of $G$.

\begin{proposition}
\label{propo:sp-radius}
Let $G$ be a distance-regular graph with adjacency and Laplacian matrices $\A$ and $\L$. Let $F_2=F_2(G)$ be its 2-token graph with adjacency and Laplacian matrices $\A(F_2)$ and $\L(F_2)$
with respective spectral radii $\rho_A(F_2)$ and $\rho_L(F_2)$.
Let $\A(F_2/\pi)$ and $\L(F_2/\pi)$ be  the quotient matrices in \eqref{A-quotient2} and \eqref{L-quotient} with respective spectral radii $\rho_A(F_2/\pi)$ and $\rho_L(F_2/\pi)$. Then, the following holds:
\begin{itemize}
\item[$(a)$]
$\rho_A(F_2)=\rho_A(F_2/\pi)$.
\item[$(b)$]
$\rho_L(F_2)\ge \rho_L(F_2/\pi)$, with equality if $G$ is bipartite.
\end{itemize}
\end{proposition}
\begin{proof}
$(a)$ This follows from Lemma \ref{regp-drg} and the known fact that if $G$ is a graph with regular partition $\pi$, then the spectral radius of $G$ is equal to the spectral radius of $G/\pi$. (The reason is that the Perron vector of $G/\pi$ lifts to the Perron vector of $G$.)

$(b)$ In this case, we know that the eigenvector $\vecv$ of the spectral radius $\rho_L(F_2/\pi)$ satisfies $\vecv^{\top} \1=0$. Then, we only can conclude that $\rho_L(F_2)\ge \rho_L(F_2/\pi)$. If $G$ is bipartite, all the closed walks rooted at a vertex $u$ have even length. Then, the diagonal entries of $\L(F_2)^{\ell}$ are the same as those of the $\ell$-th power of the signless Laplacian matrix. So, as in case $(a)$, we can apply the Perron-Frobenius Theorem, obtaining the equality between the Laplacian spectral radius of $F_2$ and that of its quotient.
\end{proof}

In fact, the eigenvalues of $\A(F_2/\pi)$ are 2 times the zeros of the so-called {\em conjugate polynomial $\overline{p}_d$ of the distance polynomial $p_d$ of $G$ (with $p_d(\A)=\A_d$, where $\A_d$ is the $d$-distance matrix of $G$). The conjugate polynomials were introduced by Fiol and Garriga in \cite{fg97}, and are defined  on the  mesh $\{\theta_0,\theta_1,\ldots,\theta_d\}$ in terms of the distance polynomials $p_0,p_1,\ldots,p_d$ as
$$
\overline{p}_i(\theta_i)=\frac{p_{d-i}(\theta_i)}{p_d(\theta_i)}\quad\mbox{for}\quad i=0,1,\ldots,d.
$$
Thus, in particular $\overline{p}_d(\theta_i)=p_d(\theta_i)^{-1}$ and, up to a constant,  $\overline{p}_d$(x) equals the characteristic polynomial of $\frac{1}{2}\A(F_2/\pi)$ (for more details, see 
C\'amara,  F\`abrega,  Fiol, and  Garriga \cite{cffg09}).
}


Let us show an example. 
\begin{example}
The Heawood graph $H$ (which is the point/line incidence graph of the Fano plane) is a bipartite distance-regular graph with $n=14$ vertices, diameter three, and intersection array $\{b_0,b_1,b_2; c_1,c_2,c_3\}=\{3,2,2; 1,1,3\}$.
The Laplacian spectral radius of $H$ is $\rho_L(H)=6$, and the algebraic connectivity of $\overline{H}$ is $\alpha_1(\overline{H})=n-\rho_L(H)=8$.
By Proposition \ref{propo:sp-radius}, the $2$-token
graph $F_2=F_2(H)$ has a regular partition $\pi$ with quotient matrix 
$$
\A(F_2/\pi)  =2\left(
\begin{array}{ccc}
0 & 1 & 0 \\
2 & 0 & 3 \\
0 & 2 & 0
\end{array}
\right),\label{A-quotient-H}
$$
and  quotient Laplacian matrix
$$
\L(F_2/\pi) =2\left(
\begin{array}{rrr}
1  & -1 & 0 \\
-2 & 5 & -3\\
0 & -2 & 2
\end{array}
\right). \label{L-quotient-H}
$$
The eigenvalues of $\A(F_2/\pi)$ are $0,\pm 4\sqrt{2}$, whereas those of $\L(F_2/\pi)$ are $0,8\pm 2\sqrt{7}$. Thus, we conclude that $\rho_A(F_2(H))=4\sqrt{2}$ 
and $\rho_{2}(H):=\rho_L(F_2(H))= 8+2\sqrt{7}$. From this and Lemma \ref{lem:alpha+rho}, we have that $\alpha_2(\overline{H})=2(n-1)-\rho_{2}(H)=18-2\sqrt{7}$, which is greater than  $\alpha_1(\overline{H})=8$. Since the algebraic connectivity of $F_2(\overline{H})$ also is 8, we get $\alpha_1(F_2(\overline{H}))=\alpha_1(\overline{H})$, as expected.\\

It is well-known that the distance polynomials satisfy the three-term recurrence
$xp_i=b_{i-1}p_{i-1}+a_i p_i+c_{i+1}p_{i+1}$ with $p_0=1$ and $p_1=x$. Then, the 
$3$-distance polynomial of $H$ is $p_3(x)=\frac{1}{3}(x^3-5x)$, so that the conjugate polynomial $\overline{p}_3$ must satisfy $\overline{p}_3(\pm 3)=p_3(\pm 3)^{-1}=\pm\frac{1}{4}$,  and $\overline{p}_3(\pm\sqrt{2})=p_3(\pm\sqrt{2})^{-1}=\mp\frac{\sqrt{2}}{2}$. This results into $\overline{p}_3(x)=\frac{1}{12}x^3-\frac{2}{3}x$, with roots $0,\pm 2\sqrt{2}$, which correspond to the eigenvalues of $\frac{1}{2}\A(F_2/\pi)$, as predicted.
\end{example}

\begin{figure}[t]
\begin{center}
\includegraphics[width=12cm]{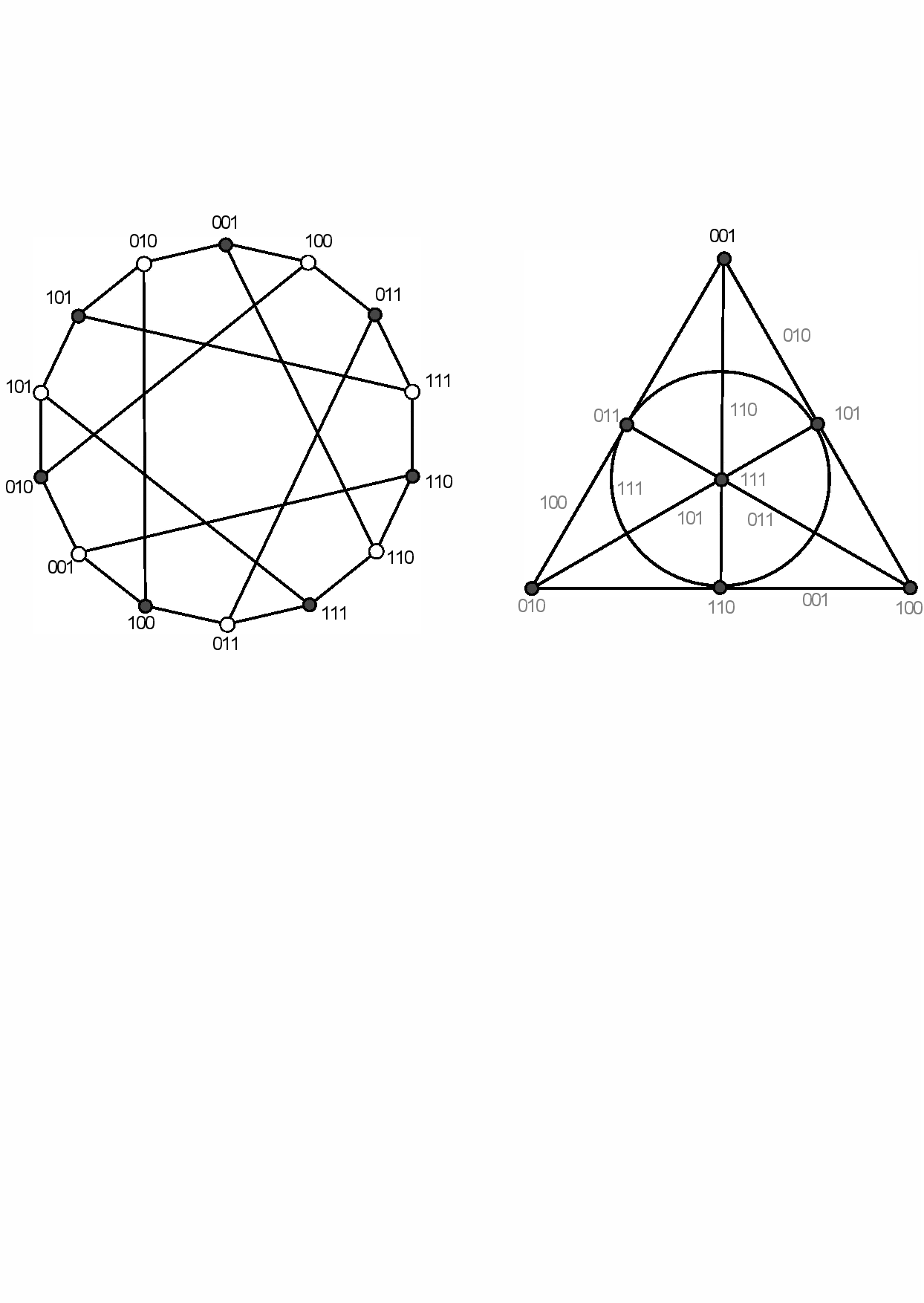}
\caption{The Heawood graph is the point-line incidence graph of the Fano plane.} 
\label{fig2}
\end{center}
\end{figure}

Other consequences of Lemma \ref{regp-drg} and Proposition \ref{propo:sp-radius} are the following.
First, from Theorem \ref{thm:Godsil},
we got the following result.
\begin{corollary}
Let ${\cal F}$ be the family of
all distance-regular graphs with diameter $d$ and the same parameters (or intersection array).
Then, every graph $G\in{\cal F}$ has 2-token graph $F_2=F_2(G)$ with the $d$  (adjacency or Laplacian) eigenvalues of $\A(F_2/\pi)$ or $\L(F_2/\pi)$ given in \eqref{A-quotient2}--\eqref{L-quotient}). In particular,  $F_2$ has spectral radii   $\rho_A(F_2)=\rho_A(F_2/\pi)$, and $\rho_L(F_2)=\rho_L(F_2/\pi)$ if it is bipartite (Proposition \ref{propo:sp-radius}).
\end{corollary}
Thus, the natural question is if, as in the case of strongly regular graphs (see Godsil \cite{g93,g95}), all distance-regular graphs with the same parameters are also cospectral (with respect to the adjacency or Laplacian matrix). In fact, notice that since every graph $G\in{\cal F}$ has the same eigenvalues, the spectrum of its $k$-token graph $F_k(G)$ also has such eigenvalues (Theorem \ref{th:sp-token}).

Moreover, from Theorem \ref{th:lew} and the interlacing theorem (see Haemers  \cite{h95}), we get the following consequence.
\begin{corollary}
Let $G$ be a distance-regular graph with (adjacency) eigenvalues $\theta_0>\theta_1>\cdots>\theta_d$. Then, the $2$-token graph $F_2(G)$ has some eigenvalues $\mu_0>\mu_1>\cdots>\mu_{d-1}$ satisfying
$$
2\theta_{i+1}\leq \mu_i \leq 2\theta_i,\qquad i=0,\ldots,d-1.
$$
\end{corollary}
\begin{proof}
Apart from the factor $2$, the matrix $\A(F_2/\pi)$  in   \eqref{A-quotient2} is a principal $d\times d$ submatrix of the intersection matrix $\B$ in  \eqref{A-quotient}. Then, the result follows by using interlacing (see Haemers \cite{h95}).
\end{proof}

\subsection{Strongly regular graphs}

Let $G$ be a (connected) strongly regular graph on $n$ vertices, which is a distance-regular graph with diameter $2$. Let $G$ have parameters $(n,d,a,c)$, that is, $G$ is $d$-regular (with $b_0=d$), $a_1=a$, and $c_2=c$. Then, its intersection matrix is
$$
\B =
\left(
\begin{array}{ccc}
0 & 1 & 0 \\
d & a & c \\
0 & d-a-1 & d-c
\end{array}
\right).\label{A-quotient-H2}
$$
Then, the $2$-token
graph $F_2=F_2(G)$ has a regular partition $\pi$ with quotient matrix 
$$
\A(F_2/\pi)  =\left(
\begin{array}{cc}
2a & 2c \\
2d-2a-2 & 2d-2c \\
\end{array}
\right).\label{A-quotient-srg}
$$
Such a regular partition was given by 
Audenaert, Godsil,  Royle, and  Rudolph in \cite{agrr07}, and noted that the adjacency eigenvalues of $\A(F_2/\pi)$ are 
$$
\theta_{1,2}=d+(a-c)\pm \sqrt{[d-(a-c)]^2-4c}.
$$
They also commented that the positive eigenvalue $\theta_1$ has a positive eigenvector (Perron vector) and, so, it corresponds to the (adjacency) spectral  radius $\rho_A(F_2(G))$.

In contrast, the quotient Laplacian matrix is 
$$
\L(F_2/\pi) =\left(
\begin{array}{cc}
2c  & -2c \\
-2d+2a +2& 2d-2a-2
\end{array}
\right), \label{L-quotient-srg}
$$
which has eigenvalues $\lambda_1=0$ and $\lambda_2=2(d-1)-2(a-c)$. Now, the eigenvector of $\lambda_2$ is orthogonal to $\1$. Then, we can only conclude that the Laplacian  spectral radius of $F_2(G)$ satisfies
\begin{equation}
\label{bound-rho}
\rho_L(F_2(G))\ge 2(d-1)-2(a-c).
\end{equation}
For instance, the cycle on five vertices $C_5$ is strongly regular with parameters $(5,2,0,1)$. Its Laplacian spectral radius is (approximately)
$\rho_L(F_2(G))=6.2361$, whereas the lower bound in \eqref{bound-rho} gives $4$.

\section*{Statements and Declarations}

The authors have no competing interests.


\end{document}